\documentclass[a4paper, twoside, 11pt]{amsart}

\usepackage[T1]{fontenc}
\usepackage[utf8]{inputenc}
\usepackage[english]{babel}
\usepackage[all]{xy}

\usepackage{lmodern}
\usepackage{amssymb,amsmath,amsthm,amsfonts,mathtools}
\usepackage[hidelinks]{hyperref}
\usepackage[nobysame,msc-links]{amsrefs}

\RequirePackage{mathrsfs}
\usepackage{enumitem}

\usepackage{booktabs}
\usepackage[font=small, margin=40pt,figureposition=below]{caption}

\usepackage{comment}

\usepackage{tikz-cd}

\newtheorem{theorem}{Theorem}[section]
\newtheorem*{theorem*}{Theorem}
\newtheorem{corollary}[theorem]{Corollary}
\newtheorem{lemma}[theorem]{Lemma}
\newtheorem{proposition}[theorem]{Proposition}
\theoremstyle{definition}
\newtheorem{definition}[theorem]{Definition}
\newtheorem{example}[theorem]{Example}
\newtheorem{remark}[theorem]{Remark}
\newtheorem{remarks}[theorem]{Remarks}

\setlist[enumerate, 1]{label = (\roman*), ref = \roman*}

\newcommand{\BA}{{\mathbb {A}}}

\newcommand{\BC}{{\mathbb {C}}}
\newcommand{\BD}{{\mathbb {D}}}

\newcommand{\BN}{{\mathbb {N}}}

\newcommand{\BP}{{\mathbb {P}}}
\newcommand{\BQ}{{\mathbb {Q}}}
\newcommand{\BR}{{\mathbb {R}}}

\newcommand{\BZ}{{\mathbb {Z}}}

\newcommand{\CA}{{\mathcal {A}}}

\newcommand{\CL}{{\mathcal {L}}}

\newcommand{\CO}{{\mathcal {O}}}

\newcommand{\CU}{{\mathcal {U}}}

\newcommand{\CW}{{\mathcal {W}}}

\newcommand{\CZ}{{\mathcal {Z}}}

\newcommand{\RH}{{\mathrm {H}}}

\newcommand{\genus}{{\rm genus}}
\newcommand{\sing}{{\rm Sing}}
\newcommand{\vol}{{\rm vol}}
\newcommand{\supp}{{\rm supp\;}}

\newcommand{\ep}{\varepsilon}
\newcommand{\la}{\lambda}
\newcommand{\La}{\Lambda}
\newcommand{\Ma}{\mathbf{M}}

\DeclareMathOperator{\NS}{NS}
\DeclareMathOperator{\Aut}{Aut}

\usepackage{color}

\begin{document}

	\title[Loxodromic K3 surfaces and complex dynamics]{Rational points on loxodromic K3 surfaces and complex dynamics
	}
	
	\author{Damián Gvirtz-Chen}
	\address{School of Mathematics and Statistics,
		University of Glasgow,
		University Place,
		Glasgow G12 8QQ,
		United Kingdom}
	\email{damian.gvirtz@glasgow.ac.uk}
	
	\author{Zhizhong Huang}
	\address{State Key Laboratory of Mathematical Sciences, Academy of Mathematics and Systems Science, Chinese Academy of Sciences, Beijing 100190, China}
	\email{zhizhong.huang@yahoo.com}

	\author{Zhuchao Ji}
	\address{Institute for Theoretical Sciences, Westlake University, Hangzhou 310030, China}
	\email{jizhuchao@westlake.edu.cn}
	
	\author{Junyi Xie}
	\address{Beijing International Center for Mathematical Research, Peking University, Beijing 100871, China}
	\email{xiejunyi@bicmr.pku.edu.cn}

	\date{\today}
	
	\begin{abstract}
		We establish an analogue of the Hilbert irreducibility theorem for all K3 surfaces over number fields that admit an automorphism of positive entropy. Our method draws ingredients from both algebraic geometry and complex dynamics.
	\end{abstract}
	
	\maketitle
	
	\tableofcontents
	
	\section{Introduction}
	\subsection{Hilbert property and K3 surfaces}
	Let $X$ be a geometrically integral variety over a field $k$.  We say that the set $X(k)$ is \emph{thin} (cf. \cite[Definition 3.1.1]{Serre}) if there exists a generically finite $k$-morphism $f:Y\to X$ where $Y$ is a $k$-variety (which may be reducible), such that $X(k)\subseteq f(Y(k))$ and $f$ admits no $k$-rational section. The last condition means that the degree of $f$ is at least 2 when restricted to any irreducible component of $Y$ of the same dimension as $Y$ that is mapped dominantly to $X$. The following density property of rational points has been originally motivated by the inverse Galois problem (cf. \cite{Serre}). 
	\begin{definition}
		We say that $X$ has
		\begin{itemize}
			\item (cf. \cite[p. 189]{CT-Sansuc}, \cite[Definition 3.1.2]{Serre}) the \emph{Hilbert property} if the set $X(k)$ is not thin;
			\item (cf. \cite[Definition 1.1]{Gvirtz-Mezzedimi}) the \emph{potential Hilbert property} if there exists a finite extension $K/k$ such that the $K$-variety $X_K$ has the Hilbert property.
		\end{itemize}
	\end{definition}

	If the base field $k$ is a number field, it follows from the Hilbert irreducibility theorem \cite{Hilbert} (see e.g. \cite[Theorem~3.4.3]{Serre}) that projective spaces have the Hilbert property. The Hilbert property (resp. the potential Hilbert property) is strictly stronger than the set of rational points being Zariski dense (resp. potentially Zariski dense).  
	The basic (counter-) examples are abelian varieties, as a consequence of the weak Mordell--Weil theorem. It has been conjectured by Corvaja and Zannier \cite[Question-Conjecture~1]{Corvaja-Zannier} that these two properties are equivalent for algebraically simply connected varieties, in particular K3 surfaces. K3 surfaces are also important classes of \emph{special varieties} introduced by Campana \cite[Definition 2.41]{Campana}, a notion which fits into his broader program of characterizing potential density of rational points (cf. \cite[Conjecture III$_A$]{Campana} and \cite[Conjecture 1.3]{Gvirtz-Mezzedimi}).
    
	The Hilbert property clearly implies the potential Hilbert property. Although the reverse direction, assuming Zariski density of rational points, is not yet known, both properties are stabilized by finite field extensions (cf. \cite[Proposition 3.2.1]{Serre}), and, as suggested by the work \cite{Gvirtz-Mezzedimi}, the potential Hilbert property allows more geometric properties of the variety to be touched.
    
    Previous approaches towards the (potential) Hilbert property for K3 surfaces are through either the existence of multiple elliptic fibrations \cites{Corvaja-Zannier,Demeio,Gvirtz-Mezzedimi,GCM24,DGC25}, or combining the weak Hilbert property for abelian varieties \cite{CDJLZ}  with pulling back a distinguished family of rational curves for Kummer surfaces associated to the Jacobians of genus two curves \cite{Gvirtz-Huang}. The goal of this article is to make use of dynamical properties of the automorphisms on the surface, and bring in tools from complex dynamics into the study of the (potential) Hilbert property.
	
	\medskip
	\subsection{Statement of the  main result}
	
    To state our result, let us first recall some terminology from dynamical systems (see e.g. \cite{CantatSurvey}). Let $X$ be a projective surface over $\BC$. The automorphism group $\operatorname{Aut}(X)$ acts on $H^*(X,\BC)$. For $f\in\operatorname{Aut}(X)$, recall that the \emph{dynamical degree} $\lambda_{f}$ is the spectral radius of $f^*$ acting on $H^{1,1}(X,\BR)$. We say that $f$ has \emph{positive entropy}  if $\lambda_{f}>1$.  It follows from \cite[\S2.2.2 \S2.3.2]{CantatSurvey} that $\lambda_{f}$ is the unique eigenvalue of $f^*$ acting on $H^2(X,\BC)$ with absolute value $>1$, which is a positive real algebraic number. 
    
Now let $X$ be a projective surface over a number field $k$. Via an embedding $\overline{k}\hookrightarrow \BC$,  $\operatorname{Aut}(X_{\overline{k}})$ embeds as a subgroup of $\operatorname{Aut}(X_{\BC})$. We say that an element in $\operatorname{Aut}(X_{\overline{k}})$  has \emph{positive entropy} if it has positive entropy in $\operatorname{Aut}(X_{\BC})$.  We say $X$ is  \emph{loxodromic} if there exists an element in  $\operatorname{Aut}(X_{\overline{k}})$ having positive entropy.

    A \emph{K3 surface} is a smooth, projective, geometrically integral surface $X$, with trivial canonical divisor, and $H^1(X,\mathcal{O}_X) = 0$. 
    Yu \cite[Theorem 1.2]{Yu} gives a lattice-theoretic characterization of complex projective K3 surfaces admitting an automorphism of positive entropy.

\medskip
	
	Our main result is the following.
	\begin{theorem}\label{thm:main}
		Let $X$ be a loxodromic K3 surface over a number field $k$. Then $X$ has the potential Hilbert property.
	\end{theorem}
	\begin{remarks}\label{rmk:tomainthm}
        \hfill
        \begin{enumerate}
		\item A classical result of Bogomolov--Tschinkel \cite[Theorem 4.10]{Bogomolov-Tschinkel} shows that rational points on K3 surfaces with an infinite automorphism group are potentially Zariski dense. The automorphism group of any loxodromic K3 surface is clearly infinite. Our result therefore goes beyond the potential Zariski density of rational points for such surfaces.  
        
        \item      
        That any K3 surface possessing at least two elliptic fibrations geometrically has the potential Hilbert property follows from the works of Demeio \cite{Demeio} and  Gvirtz-Chen--Mezzedimi \cite{Gvirtz-Mezzedimi,GCM24}. On the other hand, if $X_{\overline{k}}$ admits only one elliptic fibration, then every automorphism has to preserve it (i.e., maps fibers of the elliptic fibration to fibers), and therefore does not have positive entropy (cf. \cite[Corollaire 2.2]{cantat01}). So Theorem \ref{thm:main} covers the new cases of (loxodromic) K3 surfaces without any elliptic fibration. In fact, our results prove the Hilbert Property over the ground field for such surfaces, whenever an automorphism $f$ of positive entropy and a non-$f$-periodic rational curve are defined over $k$. See \S\ref{se:dynamicalcriterion}  for precise statements and \S\ref{se:loxoK3woellipfi} for examples. 
       \end{enumerate}
	\end{remarks}

\medskip

\subsection{Notation and conventions}
For a complex projective surface $X$, let \[\operatorname{NS}(X):=H^{1,1}(X,\BR)\cap H^2(X,\BZ)\] be the N\'eron--Severi group. Let \[\rho(X):=\operatorname{rank}\operatorname{NS}(X)\] be the Picard rank. We use $\langle\cdot,\cdot\rangle$ to denote the intersection pairing between curves, divisors, and cohomology classes.

Let $C_1, C_2$ be two irreducible algebraic curves on a projective surface which intersect properly.  For any  $x\in C_1\cap C_2$, $\langle C_1,C_2\rangle_{x}$ denotes the local intersection number at $x$.  We say that  $C_1$ and $C_2$ \emph{intersect transversally} (cf. \cite[p. 357]{Hartshorne}) at a point $x\in C_1\cap C_2$ if the local equations of $C_1, C_2$ generate the maximal ideal of the local ring at $x$. This implies that both $C_1,C_2$ are smooth at $x$, and $\langle C_1,C_2\rangle_{x}=1$. 
Let $C_1\pitchfork C_2$ denote the subset of $C_1\cap C_2$ where the intersection is transversal.

We shall call a proper dominant generically finite morphism $\pi: Y\to X$ of degree at least two between geometrically integral smooth projective varieties $X, Y$ a \emph{cover}. For such a cover, let $R_{\pi}:=K_Y-\pi^*(K_X)$ be the \emph{ramification divisor} of the cover $\pi$ in $Y$, where $K_X$, $K_Y$ are the canonical divisors of $X$ and $Y$ respectively. Consider the Stein factorization (cf. \cite[III. Corollary 11.5]{Hartshorne}) $\pi=\pi'\circ \iota$, where $\iota:Y\to Y'$ is proper and has connected fibers, and $\pi':Y'\to X$ is finite. We denote by $D_\pi$ the \emph{branch locus} of $\pi'$ (i.e. the locus in $X$ where $\pi'$ is not \'etale). If $X$ has trivial \'etale fundamental group (in particular if $X$ is a projective K3 surface), then $D_\pi$ is always non-empty, and by the Zariski--Nagata purity \cite[Tag OBMB]{Stacks-project-Lemma-58.21.4} (note that $Y'$ is normal),  $D_\pi$ has pure of codimension one in $X$. 

Projection formulas (cf. \cite[Proposition~2.3(c), p.~34]{Fulton1998}, \cite[II.~Exercise~5.1(d)]{Hartshorne}) of the following form are used throughout. Let $\pi:Y\to X$ be a cover between smooth projective surfaces. Let $A$ be a curve (resp. divisor, cohomology class) on $X$ and let $B$ be a curve (resp. divisor, cohomology class) on $Y$. Then the intersection pairing satisfies
$$\langle\pi^*(A),B\rangle=\langle A,\pi_*(B)\rangle.$$
Now assume that $A,B$ are curves, and that $B$ is not vertical with respect to $\pi$\footnote{Here and after, a curve $B$ is called  {\em vertical} if $\pi(B)$ is a point.}.  Localizing at a point $x\in X$ where the intersection of $A$ with $\pi(B)$ is proper, we obtain the following equality of local intersection numbers
$$\sum_{y\in \pi^{-1}(x)\cap B}\langle\pi^*(A),B\rangle_y=\langle A,\pi_*(B)\rangle_x.$$
\medskip
\subsection{Strategy of the proof}
	
	The idea of our proof is summarized below.  Given an automorphism $f\in \operatorname{Aut}(X_{\overline{k}})$ with positive entropy on a projective K3 surface $X$ over a number field $k$, we can find a rational curve $C$ with an infinite $f$-orbit, upon passing to a finite extension of $k$. To establish the potential Hilbert property, inspired by the strategy of the work of Gvirtz-Chen--Huang \cite{Gvirtz-Huang}, we intend to show that for any fixed cover $\pi: Y\to X$, for all $n$ sufficiently large, $f^n(C)$, which is a rational curve, cannot be birationally lifted to a rational curve in $Y$. As the union of $f^n(C)$ is Zariski dense in $X$, we then make use of the Hilbert property for the projective line to conclude. 
    For simplicity, we explain our strategy under the generic setting that the branch locus $D_\pi\subseteq X$ is not $f$-periodic.
	\medskip
	
	{\bf Step 1: Producing transversal intersections using complex dynamics.}  The dynamical properties of $f$ can ``stretch'' both $C$ and a non-periodic irreducible curve $D\subseteq D_\pi$ to laminar currents that intersect each other transversally, cf.\ Proposition \ref{prop:trans}.  In this step, we use techniques from complex dynamics, in particular the theory of laminar currents, and the properties of Green currents and Green measures for positive entropy automorphisms. 
	
	\medskip

{\bf Step 2: Producing non-transversal intersections using algebraic geometry.} Assume by contradiction that there exists an infinite subfamily $I\subseteq \BZ_{\geq 0}$ such that $f^n(C)$ (which is a rational curve) can be birationally lifted to a rational curve in $Y$ for all $n\in I$.  We can then produce non-transversal intersections between $(f^n(C))_{n\in I}$ and a non-periodic irreducible curve $D\subseteq D_\pi$ due to ramification, cf. Proposition \ref{prop:non-trans} and Lemma \ref{le:non-trans}.  In this step, we use techniques from algebraic geometry, in particular the Hodge Index Theorem.  In view of Step 1, we will get a contradiction, which will finish the proof.

\medskip

It is worth emphasizing that transversality and non-transversality of the intersections that we obtain here are in a statistical sense: in Step 1, we prove that $100\%$ of the intersections are transversal as $n\to\infty$, while in Step 2, a positive proportion of the intersections are non-transversal as $n$ varies in $I$. This is sufficient to derive a contradiction.

\begin{remark}
    It would be interesting to deduce the (potential) Hilbert property for K3 surfaces with an infinite automorphism group but without elements of positive entropy, which would go beyond the techniques in our paper. (See \cite[Theorem 1.1]{Yu} for a classification of K3 surfaces admitting an infinite automorphism group preserving a unique elliptic fibration.) The classification of K3 surfaces having a finite automorphism group goes back to Nikulin \cite{Nikulinautomorphism} (see also \cite[Chapter~\S15.2.4]{Huybrechts}). It appears to be a challenge to prove the potential Zariski density of rational points on such surfaces. 
\end{remark}
    
    \medskip
    \subsection{Organization of the paper}
In \S\ref{se:dynamicalcriterion} we reduce the proof of Theorem \ref{thm:main} to an explicit dynamical and algebraic geometry type statement -- Theorem \ref{thm:maingroundfield} -- about pulling back the orbit of a non-periodic curve along any given cover. We show how Propositions \ref{prop:trans} and \ref{prop:non-trans} together yield Theorem \ref{thm:maingroundfield}.
In \S\ref{sec:laminar} we introduce the terminology of laminar currents and their variations, and prove some basic properties. Equipped with these, in \S\ref{sec:trans} we prove Proposition \ref{prop:trans}.
In \S\ref{sec:non-trans} we prove Proposition \ref{prop:non-trans} relying on ramification.
    \medskip
    
	\section{A criterion for liftability of orbits of curves along covers}\label{se:dynamicalcriterion} 
	
Let $X$ be a complex K3 surface.  
For $f\in \operatorname{Aut}\left(X\right)$, let $B_f\subseteq X$ denote the set of all $f$-periodic curves. If $f$ has positive entropy, then by \cite[Proposition 3.1 (1)]{Kawaguchi}, $B_f$ is the union of finitely many curves.

As the technical core for completing the proof of Theorem \ref{thm:main}, we shall establish the following:
	\begin{theorem}\label{thm:maingroundfield}
		Let $X$ be a complex K3 surface. 
        Let $f\in \operatorname{Aut}\left(X\right)$ with positive entropy. 
        Let $C\subseteq X$ be an irreducible algebraic curve
        with infinite orbit under $f$ and let $C_n:=f^n(C)$, $n\geq 0$. Let $\pi: Y\to X$ be a cover such that $D_{\pi}\not\subseteq B_f$.
        Then for $n\gg 1$, every non-vertical irreducible component of $\pi^{-1}(C_n)$ maps to $C_n$ with degree $\geq 2$. 
	\end{theorem}
	
Before proceeding to the proof of the more general Theorem \ref{thm:main}, we first deduce the following relatively simpler dynamical criterion for the Hilbert property (over the ground field).

\begin{theorem}\label{thm:hilbbfempty}
    Let $X$ be a  K3 surface over a number field $k$. Assume that there exist a rational curve $C\subseteq X$ over $k$ and $f\in \operatorname{Aut}\left(X\right)$ over $k$ with positive entropy such that $B_f=\varnothing$. Then $X$ has the Hilbert property (over $k$).
\end{theorem}
    \begin{proof}
        Granting Theorem \ref{thm:maingroundfield}, we directly apply the criterion \cite[Proposition 5.1]{Gvirtz-Huang} to conclude.
    \end{proof}

    Dealing with automorphisms with non-empty periodic loci requires a more delicate analysis. 
\medskip

    \subsection{Considerations for periodic curves}
    Let $X$ be a complex K3 surface and $G\leq\Aut(X)$ be a subgroup containing at least one element of positive entropy. Write $R_G\subseteq \NS(X)$ for the sublattice spanned by the (classes of) curves on $X$ which have finite orbit under $G$, and $B_G\subsetneq X$ for their union.

		\smallskip
		
		\begin{lemma}\label{le:simplyconnectedness}
		     Assume that $\operatorname{rank}R_G< 8$. Then the open subset $X\setminus B_G$  is simply connected.
		\end{lemma}
\begin{proof}
 We first note that the lattice $R_G$ is negative
definite: its classes are
orthogonal to the eigenclasses $\theta_f^\pm$ of a positive entropy automorphism $f$   satisfying
$f^*\theta_f^\pm=\lambda_f^{\pm1}\theta_f^\pm$ and
$(\theta_f^+ +\theta_f^-)^2>0$, so the Hodge index theorem applies 
(see \cite[Theorem~2.1 and \S2.3.2]{CantatSurvey}).

We may assume $B_G\ne\varnothing$.
By adjunction, every component of $B_G$ is a $(-2)$-curve.
There are finitely many such components, since $R_G$ has finitely many vectors of
square $-2$, and every such class is represented by at most one
irreducible curve. Their classes are linearly independent: otherwise
there would be nonzero effective divisors $A,B$ with the same class in $R_G$ but with no common
components, contradicting
$A^2=\langle A,B\rangle\geq 0$. Thus $B_G$ is an ADE configuration, and $R_G$ with the intersection pairing, is a \emph{root lattice}:
a free abelian group with a negative-definite integral symmetric
bilinear form, generated by vectors of square $-2$, called
\emph{roots}.

By \cite[Propositions~4.1 and~4.6]{Gvirtz-Mezzedimi}, any negative-definite root lattice
$R\subseteq\operatorname{NS}(S)$ generated by $(-2)$-curves on a
K3 surface $S$ and satisfying $\operatorname{rank}R<8$
is primitive in $\operatorname{NS}(S)$. Here, \emph{primitivity} means that $\NS(S)/R$ is torsion-free.
Thus, in our case, $R_G$ is primitive in $\NS(X)$.
Since $\NS(X)$ is primitive in $H^2(X,\BZ)$, $R_G$ is also primitive in $H^2(X,\BZ)$.

Let $A_{R_G}:=R_G^\vee/R_G$ be the \emph{discriminant group}, where
$$R_G^\vee:=\{v\in R_G\otimes\BQ:v\cdot R_G\subset\BZ\}$$ is  the \emph{dual lattice}.  Write $\ell$
for the minimal number of generators of this finite group and $r$ for the rank of $R_G$.
We have
\[
 r\leq18,\qquad
 \ell(A_{R_G})\leq r=\min\{r,20-r\}.
\]
We then conclude by  \cite[Theorem~4.3]{Shimada-Zhang} that $\pi_1(X\setminus B_G)=1$.
\end{proof}

\medskip

\begin{proposition}\label{prop:hilboverk}
 Let $X$ be a K3 surface over a number field $k$ and let  $G\leq\Aut_k(X)$ be a subgroup containing at least one element of positive entropy. Assume that $X_\BC\setminus B_G$ is simply connected and that there exists a rational curve $C$ on $X$ over $k$ such that $C_\BC\nsubseteq B_G$. Then $X$ has the Hilbert Property (over $k$).
\end{proposition}

\begin{proof}
Choose an element $f\in G$ of positive entropy. Since the $G$-orbit
of $C$ is infinite and $B_f$ has only finitely many irreducible
components, there exists $h_0\in G$ such that
\[
h_0(C)\not\subseteq B_f.
\]
Replacing $C$ by $h_0(C)$ and keeping $f$ unchanged, we may assume
that $C$ is not $f$-periodic.

 Now given any collection of finitely many covers $(\pi_i:Y_i\to X)_{1\leq i\leq r}$ over $k$, we intend to show that there exist $h\in G$ and $n_0\in \BN$, such that the following holds: Let $C'_n:=h(C_n)$. Then for 
     every $n\geq n_0$ and for every $1\leq i\leq r$,  every non-vertical  irreducible component of $\pi_i^{-1}(C'_n)$ is mapped to $C'_n$ with degree $\geq 2$.  
     From this, by using the  Hilbert property of $\mathbb{P}_k^1$, we deduce the  Hilbert property for $X$ by appealing to \cite[Proposition 5.1]{Gvirtz-Huang}.

In what follows we base change everything to $\BC$ (and omit the subscript). Since $X\setminus B_G$ is simply connected,
we have $D_{\pi_i}\not\subset B_G$ for every $1\leq i\leq r$.
Choose an irreducible component $D_i\subset D_{\pi_i}$ such that $D_i\not\subset B_G.$ By the definition of $B_G$, the $G$-orbit of each $D_i$ is infinite. 

We claim that there exists a single $h\in G$ such that $D_i\not\subset h(B_f)$ for all $1\le i\le r$. Write $B_f=\cup_{j=1}^{s} B_j$,
where $B_j$ are the irreducible components of $B_f$.
For each pair $(i,j)$, define
\[
T_{ij}:=\{h\in G:h(B_j)=D_i\}.
\]
If $T_{ij}$ is nonempty, we fix $h_{ij}\in T_{ij}$. Then
\[
T_{ij}=\operatorname{Stab}_G(D_i)\,h_{ij},
\qquad
[G:\operatorname{Stab}_G(D_i)]
=\#(G\cdot D_i)=\infty.
\]
By \cite[Lemma~4.1]{Neumann1954}, a group cannot be covered by
finitely many cosets of subgroups of infinite index. Hence every \[h\in G\setminus\bigcup_{i=1}^{r}\bigcup_{j=1}^{s}T_{ij}\]  satisfies the claim.

Let
\[
g:=h\circ f\circ h^{-1},\qquad C':=h(C).
\]
Then $g$ has positive entropy, $C'$ has an infinite $g$-orbit,
and $B_g=h(B_f)$. Consequently, $D_{\pi_i}\not\subset B_g\;\text{for all }1\le i\le r.$ Let $$C'_n:=g^n(C')=h(C_n).$$ We may therefore apply Theorem~2.1 to all the covers $\pi_i$
using the same automorphism $g$ and the same sequence $(C'_n)$ of rational curves, and we conclude that there exists a common
$n_0$ such that, for every $n\ge n_0$ and every $1\leq i\leq r$, every
non-vertical irreducible component of $\pi_i^{-1}(C'_n)$
maps to $C'_n$ with degree at least two. Recall that the automorphisms $f,g,h$, and hence the sequence $(C_n')$, are all defined over the ground field. This finishes the proof.
\end{proof}

	\subsection{Proof of Theorem \ref{thm:main}}
    The intersection pairing on $\operatorname{NS}(X_{\overline{k}})$ has signature $(1,\rho(X_{\overline{k}})-1)$. By \cite{PyatetskiiShapiroShafarevich} (cf. also \cite[Proposition 11.1.3]{Huybrechts}), elliptic fibrations on $X_{\overline{k}}$ are in one-to-one correspondence with primitive nef elements $v\in\operatorname{NS}(X_{\overline{k}})$ satisfying $v^2=0$. Such an isotropic element always exists if $\rho(X_{\overline{k}})\geq 5$.
    As explained in Remarks \ref{rmk:tomainthm} (ii), either every element of  $\operatorname{Aut}(X_{\overline{k}})$ has zero entropy, or $X_{\overline{k}}$ admits at least two elliptic fibrations, for which we appeal to the work \cite{GCM24}. We may henceforth assume  $\rho(X_{\overline{k}})<5$, which clearly implies $\operatorname{rank}R_{\Aut(X_{\overline{k}})}< 8$.

     By \cite[Proposition 2.2]{Sterk} (cf. also \cite[Corollary 15.2.4]{Huybrechts}), the group $\operatorname{Aut}\left(X_{\overline{k}}\right)$ is finitely generated. By \cite[Proof of Theorem 4.10]{Bogomolov-Tschinkel}, there exists a (possibly singular) rational curve $C$ on $X_{\overline{k}}$ such that its orbit under $\operatorname{Aut}\left(X_{\overline{k}}\right)$ is infinite. Upon passing to a finite extension of $k$, we may assume that every  automorphism on $X$ and the curve $C$  are defined over the ground field $k$. This thus confirms the potential Hilbert property for $X$ by combining Lemma \ref{le:simplyconnectedness} and Proposition \ref{prop:hilboverk}.
	\qed
	
	\medskip
	
	\subsection{Roadmap of proof of Theorem \ref{thm:maingroundfield}}
	The rest of our paper is devoted to the proof of Theorem \ref{thm:maingroundfield}. Our plan will be split into two parts. 
    
    \smallskip
    \begin{proposition}\label{prop:trans} In the setting of Theorem \ref{thm:maingroundfield}, let $D\subseteq D_\pi$ be an irreducible component such that $D\not\subseteq B_f$. Then $\langle C_n, D\rangle>0$ for $n$ large enough, and we have \[\lim_{n\to+\infty} \frac{\#(C_n\pitchfork D)}{\langle C_n, D\rangle}=1.\]
	\end{proposition}

\medskip

The proof of Proposition \ref{prop:trans} will be given in \S\ref{sec:trans}.

	\medskip
		
	\begin{proposition}\label{prop:non-trans}
	In the setting of Theorem \ref{thm:maingroundfield}, assume that there exists an infinite subfamily $(C_n)_{n\in I}$, $I\subset \BZ_{\geq0}$, such that for all $n\in I$, there exists an irreducible component $E_n\subseteq \pi^{-1}(C_n)$ such that $\pi|_{E_n}:E_n\to C_n$ is birational.
    Let $D\subseteq D_\pi$ be an irreducible component such that $D\not\subseteq B_f$. Let $R\leq R_\pi$ be the effective divisor which is the sum (with the same multiplicities) of all irreducible components $Z$ of $\operatorname{supp}R_\pi$ such that $\pi(Z)=D$. Then  we have $\langle E_n,R\rangle>0$ and $\langle C_n,D\rangle>0$ for all $n$ large enough, and
		\[\liminf_{\substack{n\to +\infty\\ n\in I}}\frac{\langle E_n,R\rangle}{\langle C_n,D\rangle}>0.\]
	\end{proposition}
	
	\medskip
	
	The proof of Proposition \ref{prop:non-trans} will be given in \S\ref{sec:non-trans}.

	\medskip
		
	We next combine  Propositions \ref{prop:trans} and \ref{prop:non-trans} to prove Theorem \ref{thm:maingroundfield}. The following basic observation is crucial for our analysis.
	\begin{lemma}\label{le:non-trans}
		In the setting of Proposition \ref{prop:non-trans}. Let $y\in E_n\cap \supp R$. Then the intersection $C_n\cap D$ at $\pi(y)$ is not transversal. \end{lemma}
	\begin{proof}
		Since $E_n$ is not vertical, by the projection formula, we have 
		\[\langle C_n,D\rangle_{\pi(y)}=\langle\pi_*(E_n), D\rangle_{\pi(y)}\geq\langle E_n, \pi^*(D)\rangle_{y}.\]
        Every summand of $R$ appears with multiplicity at least $2$ in $\pi^*(D)$.
        By the assumption  $y\in E_n\cap \supp R$,  it follows that $\langle E_n, \pi^*(D)\rangle_{y}\geq 2$, 
        hence $\langle C_n,D\rangle_{\pi(y)}\geq 2$. This implies that the intersection of $C_n$ and $D$ at $\pi(y)$ is not transversal.
	\end{proof}
    \begin{example}
        As a basic illustration of Lemma \ref{le:non-trans}, let us consider the degree $2$ endomorphism $\pi:(s,t)\mapsto (s^2,t)$ of $\BA^2$. The branch curve $D$ is $s=0$, and the curve $E:s=t$ is mapped birationally onto its image $C:s=t^2$. We see that $D$ is tangent to $C$ at $(0,0)$, although $\pi^{-1}(D)$ intersects $E$ at $(0,0)$ transversally. See the following figure.
    \begin{figure}[h]\centering\begin{minipage}{0.48\textwidth}\centering\includegraphics[scale=0.3]{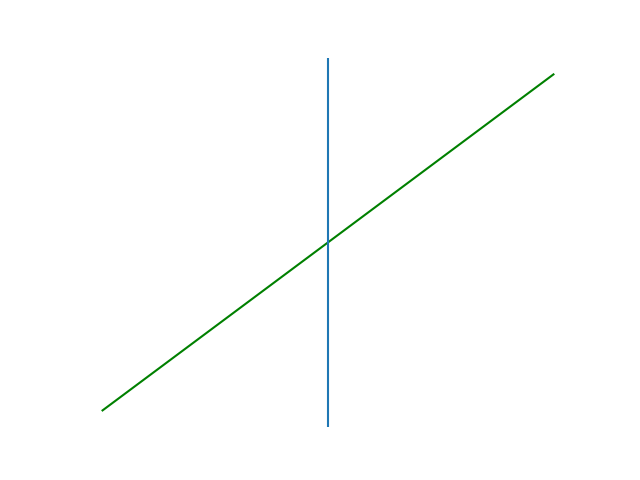}
			\end{minipage}
		\begin{minipage}{0.48\textwidth}\centering\includegraphics[scale=0.3]{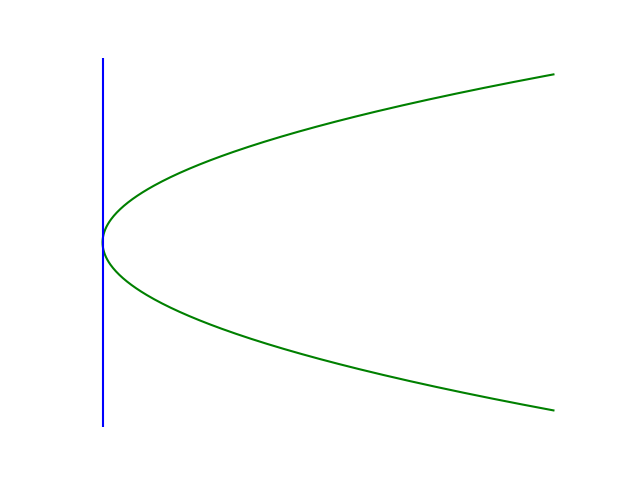}
			\end{minipage}\caption{\small The images of the curves $s=t$ and $s=0$ under the endomorphism $(s,t)\mapsto (s^2,t)$ in $\BA^2$}\label{fig:xtox2}
            \end{figure}
        
        To further illustrate the condition on $y$ in Lemma \ref{le:non-trans}, consider the degree~$3$ \emph{non-Galois} endomorphism $\pi:(s,t)\mapsto(s(s^2+st+1),t)$ of $\BA^2$. The curve $C:s=0$ is mapped isomorphically to itself under $\pi$, and $\pi(C)$ intersects the branch locus of $\pi$ transversally in two points $(0,\pm 2)$, while $C$ does not intersect the ramification locus $3s^2+2st+1=0$. 
	    \end{example}

	\begin{proof}[Proof of Theorem \ref{thm:maingroundfield}]
		We argue by contradiction. Suppose that the conclusion of Theorem \ref{thm:maingroundfield} is not true.  Then by Proposition \ref{prop:non-trans}, there exists an infinite subfamily $(C_n)_{n\in I}$ such that for all $n\in I$, $\langle C_n, D\rangle>0$, $\langle E_n, R\rangle>0$, and 
		\begin{equation}\label{eq:liminf1}
		  \liminf_{\substack{n\to +\infty\\ n\in I}}\frac{\langle E_n,R\rangle}{\langle C_n,D\rangle}>0.  \end{equation}	
		
		Now for all $n\in I$, by Lemma \ref{le:non-trans},  all the intersection points of $C_n\cap D$ that are contained in $\pi(E_n\cap \supp R)$ are  proper and non-transversal. Since $E_n$ is not vertical, we have
		\begin{equation}\label{eq:liminf2}
        \begin{split}
        \langle C_n, D\rangle - \#(C_n\pitchfork D) &\geq\sum_{x\in\pi(E_n\cap \supp R)}\langle C_n,D\rangle_x\\&\geq \sum_{y\in E_n\cap \supp R}\langle E_n,\pi^*(D)\rangle_{y}\\ &\geq \langle E_n,R\rangle.
        \end{split}
		\end{equation}

	Combining \eqref{eq:liminf1} and  \eqref{eq:liminf2}, we have 
    \begin{align*}
        \limsup_{\substack{n\to +\infty\\ n\in I}}\frac{\#(C_n\pitchfork D)}{\langle C_n, D\rangle}&\leq \limsup_{\substack{n\to +\infty\\ n\in I}}\left(1-\frac{\langle E_n,R\rangle}{\langle C_n,D\rangle}\right)\\ &=1-\liminf_{\substack{n\to +\infty\\ n\in I}}\frac{\langle E_n,R\rangle}{\langle C_n,D\rangle}\\ &<1,
    \end{align*}
	which contradicts Proposition \ref{prop:trans}. This finishes the proof. 	\end{proof}		
	
	\medskip

	\section{Laminar currents and intersection theory}\label{sec:laminar}
	The purpose of this section is to recall the theory of laminar currents and intersection theory on surfaces, which are necessary for the proof of Proposition \ref{prop:trans}.  Many results in this section stem from the theory of laminar currents   originating in the celebrated work of Bedford--Lyubich--Smillie \cite{Bedford-Lyubich-Smillie}, and its further developments due to Dujardin \cite{DujardinMathann, dujardin04, dujardin06}. (See \cite{dujardin23} for a nice introduction to the theory of laminar currents.) See also works of Cantat \cite{cantat01}, De Th\'elin \cite{dethelin04, dethelin08}, Dinh \cite{dinh05}, Diller--Dujardin--Guedj \cite{diller11},   and \cite{ji2026genus} for further developments.
	
	\subsection{Uniformly laminar currents}
Let $(X,\omega)$ be a connected Hermitian surface.  Let $\CL$ be a family of smooth connected closed analytic curves in $X$. By integrating smooth real-valued $(1,1)$-forms on $X$ over $\Gamma\in\CL$, the family $\CL$ can be seen as a subset  of the space of positive closed $(1,1)$-currents on $X$. We equip $\CL$ with the subspace topology induced by the weak topology on the space of currents. Since curves in $\CL$  are smooth and have multiplicity one, convergence in $\CL$ implies local $C^\infty$-convergence of the corresponding analytic graphs.

\begin{definition}[Uniformly laminar family]
The family $\CL$ is called a {\em uniformly laminar family} if $\CL$ is a compact metrizable space, and any two different curves in $\CL$ have empty intersection. 
\end{definition}
 Note that if $\CL$ is a compact metrizable space, then any finite Borel measure on $\CL$ is automatically a Radon measure.
\begin{example}
Let $X=\BD\times \BD$ be the unit bi-disk.  Let $(\phi_\la)_{\la\in \La}:\BD\to \BD$ be a family of holomorphic maps which is compact in the compact-open topology.    Moreover assume that for every $\la_1\neq \la_2\in \La$ and $z\in \BD$, $\phi_{\la_1}(z)\neq \phi_{\la_2}(z)$.  Let $\CL$ be the family of graphs of $\phi_\la$. Then $\CL$ is a uniformly laminar family. 
\end{example}

 \begin{definition}[Uniformly laminar currents]
 A positive closed $(1,1)$-current $S$ on $X$ is called a {\em uniformly laminar current}, if there exists a uniformly laminar family $\CL$ and a finite Borel measure $\nu$ on $\CL$ such that
 $$S=\int_{\CL} [\Gamma]\;d\nu(\Gamma),\footnote{Here and after, the notation $[\cdot]$ denotes the corresponding current of a curve.}$$
We call the above expression an \emph{integral representation} of $S$. 
 \end{definition}

Roughly speaking, a uniformly laminar current is an average of smooth curves coming from a uniformly laminar family.
\medskip

By \cite[Lemma 6.4]{Bedford-Lyubich-Smillie}, we have:
\begin{lemma}\label{lem:BLStransversal}
Let $V$ be a connected smooth closed analytic curve in  $X$, and $\CL$ be a uniformly laminar family. Assume that for a curve  $\Gamma_0\in \CL$, $\Gamma_0$ and $V$ intersect properly and non-transversally at $x\in X$. Then there exist  open neighborhood $\CU\subseteq \CL$ of $\Gamma_0$ and  open neighborhood $\Omega\subseteq X$ of $x$, such that for every $\Gamma\in \CU$, $\Gamma\neq \Gamma_0$, the intersection of $\Gamma$ and $V$ in $\Omega$ is non-empty, and all intersection points are transversal. 
\end{lemma}

As a corollary, we have:
\begin{corollary}\label{cor:transversal}
Let $V$ be a connected smooth closed analytic curve in  $X$, $\CL$ be a uniformly laminar family, and $K\subseteq X$ be a compact subset. Then
$$\#\left\{\Gamma\in \CL: \Gamma \;\text{and $V$ have a non-transversal intersection point in $K$}\right\}<+\infty.$$
In particular, the set  
$$\left\{\Gamma\in \CL: \Gamma \;\text{and $V$ have a non-transversal intersection point in $X$}\right\}$$ is countable.
\end{corollary}
\begin{proof}
Let $Z\subseteq \CL\times K$ be the subset containing all pairs $(\Gamma, x)$ such that $\Gamma$ and $V$ intersect   non-transversally at $x$.  Since the limit of non-transversal intersection points is again a non-transversal intersection point, $Z$ is closed.    By Lemma \ref{lem:BLStransversal},  $Z$ is a discrete set. Since $\CL\times K$ is compact, $Z$ is a finite set, which implies the first conclusion of Corollary \ref{cor:transversal}.  Since $X$ can be covered by countably many compact subsets, the second conclusion also holds.
\end{proof}

\subsection{Transversal geometric intersection of uniformly laminar currents}
Let $(X,\omega)$ be a connected Hermitian surface. Let $S$, $T$ be two positive closed $(1,1)$-currents on $X$.  We first recall the classical definition of the intersection $T\wedge S$ by using local potentials.   The {\em trace measure} of $S$ is by definition the measure $\sigma_S:=S\wedge \omega$.  On a local chart, we can write $T=dd^c u$, where $u$ is a plurisubharmonic function.  

\begin{definition}\label{def:admissible} 
    We say that $T\wedge S$ is {\em admissible} if on each local chart, the local potential $u$ of $T$ satisfies the integrability condition. In this case, we can globally define $T \wedge S$ via $dd^c u\wedge S:=dd^c(uS)$ on each local chart, which is a Radon measure on $X$. The \emph{total mass} is $\Ma(T\wedge S):=\int_X d(T\wedge S)$.
\end{definition} 
 If $T$ has bounded local potential, then $T\wedge S$ is always admissible. In particular, $T\wedge S$ is admissible if it has continuous bounded potential. For two closed analytic curves $\Gamma_1,\Gamma_2$ in $X$,  $[\Gamma_1]\wedge [\Gamma_2]$ is admissible if and only if $\Gamma_1$ intersects $\Gamma_2$ properly, and this happens exactly when no irreducible components of $\Gamma_1$ and $\Gamma_2$ coincide.

\medskip

\begin{definition}[Diffuse uniformly laminar currents]
	We say that  a uniformly laminar current $S=\int_{\CL} [\Gamma]\;d\mu(\Gamma)$ is {\em diffuse} if $\mu$ has no atoms, i.e. $\mu$ put no mass on every point. 
\end{definition}

If $S$ has bounded local potential, then $S$ is diffuse. 

\medskip

Let $\CL_1$, $\CL_2$ be two uniformly laminar families on $X$.    Let $$T=\int_{\CL_1} [\Gamma]\;d\mu(\Gamma) $$and $$S=\int_{\CL_2} [\Gamma']\;d\nu(\Gamma)$$ be two  uniformly laminar currents on $X$. 

\medskip

By \cite[Th\'eor\`eme 3.1]{dujardin04} or \cite[Proposition 2.11]{ji2026genus}, we have:
\begin{proposition}\label{prop:geometric intersection}
  Assume that $T\wedge S$ is admissible.  Then for $(\mu\otimes\nu)$-a.e. $(\Gamma,\Gamma')$, $\Gamma$ intersects $\Gamma'$ properly, and we have 
$$T\wedge S=\int_{\CL_1\times \CL_2} [\Gamma]\wedge [\Gamma']\;d(\mu\otimes\nu)(\Gamma,\Gamma').$$
\end{proposition}

\medskip

\begin{definition}
The \emph{transversal geometric intersection} of two uniformly laminar currents $T$ and $S$ is defined by 
$$T\dot\wedge S:=\int_{\CL_1\times \CL_2} [\Gamma\pitchfork \Gamma']\;d(\mu\otimes\nu)(\Gamma,\Gamma'),$$
  where $[\Gamma\pitchfork \Gamma']$ is the sum of Dirac masses of all transversal intersection points in $\Gamma\cap \Gamma'$. This is a positive Borel measure on $X$.
\end{definition}

We note that the measure $T\dot\wedge S$ is well defined even if the intersection $T\wedge S$ is not admissible. However, $T\dot\wedge S$ may not be Radon.  When $T\wedge S$ is admissible, we have $T\dot\wedge S\leq T\wedge S.$ The equality holds if either $T$ or $S$ is diffuse:

\begin{proposition}\label{prop:intersection equality}
Assume that $T\wedge S$ is admissible and $S$ is diffuse. Let $\CZ\subseteq \CL_1\times \CL_2$ be the subset containing pairs $(\Gamma,\Gamma')$ such that $\Gamma\cap \Gamma'$ admits a non-transversal intersection point.  Then $(\mu\otimes \nu)(\CZ)=0$.  In  particular, 
$$T\dot\wedge S=T\wedge S.$$
\end{proposition}
\begin{proof}
 Fix $\Gamma\in \CL_1$. Let $\CZ_\Gamma\subseteq \CL_2$ be the subset containing all $\Gamma'$ such that $\Gamma\cap \Gamma'$ admits a non-transversal intersection point. Thus $\CZ=\cup_{\Gamma\in \CL_1}\left\{ \Gamma\right\}\times \CZ_\Gamma$.  By Corollary \ref{cor:transversal}, for every $\Gamma\in \CL_1$, $\CZ_\Gamma$ is a countable set.  Since $S$ is diffuse, $\nu(\CZ_\Gamma)=0$.  For a compact subset $K\subseteq X$, let $\CZ_K$ be the subset of $\CZ$ containing pairs $(\Gamma,\Gamma')$ that admit a non-transversal intersection point in $K$.  Then $\CZ_K$ is closed. Since $X$ can be covered by countably many compact subsets, the set $\CZ$ is a countable union of closed subsets, hence is Borel measurable.   By Fubini's theorem, $(\mu\otimes \nu)(\CZ)=0$.  By Proposition \ref{prop:geometric intersection} we get the desired equality $T\dot\wedge S=T\wedge S.$ 
\end{proof}

\medskip

We next introduce an analogue but weaker notion of uniformly laminar family, for which we allow self-intersections. 
\begin{definition}[Uniformly woven family]
 Let $\CW$ be a family of smooth connected closed analytic curves in $X$. We call $\CW$ a {\em uniformly woven family} if $\CW$ is a compact metrizable space.
\end{definition}

Let $\CW_1$ (resp. $\CW_2$) be a uniformly woven family.   Let $\CL_1\subseteq \CW_1$ (resp. $\CL_2\subseteq \CW_2$) be  a  subfamily which is uniformly laminar.   Let $\mu$ (resp. $\nu$)  be a finite Borel measure supported on $\CL_1$ (resp. $\CL_2$).   Let $\CL_{1,n}\subseteq \CW_1$ (resp. $\CL_{2,n}\subseteq\CW_2$) be a sequence of uniformly laminar family on $X$. 
Let $\mu_n$ (resp. $\nu_n$) be a sequence of finite Borel measures supported on $\CL_{1,n}$ (resp. $\CL_{2,n}$) such that $\mu_n\to \mu$ (resp. $\nu_n\to \nu$) in the weak-* topology.  Let $$T_n:= \int_{\CW_1} [\Gamma]\;d\mu_n(\Gamma), \;\;T:= \int_{\CW_1} [\Gamma]\;d\mu(\Gamma),$$ and $$S_n:= \int_{\CW_2} [\Gamma']\;d\nu_n(\Gamma'), \;\;S:= \int_{\CW_2} [\Gamma']\;d\nu(\Gamma')$$ be the corresponding uniformly laminar currents.  Then $T_n\to T$ and $S_n\to S$ in the sense of currents. 
We have the following crucial lower estimate for the integration of a non-negative continuous function with compact support.

\begin{proposition}\label{prop:lower mass estimate}
	Notation as above. Assume that $T\wedge S$ is admissible and $S$ is diffuse. Let $\phi: X\to \BR$ be a non-negative continuous function with compact support. Then
	$$\liminf_{n\to+\infty} \int_X \phi\;d(T_n\dot\wedge S_n)\geq \int_X \phi\; d(T\dot\wedge S)= \int_X \phi\; d(T\wedge S).$$
\end{proposition}
\begin{proof}
	
	By Proposition \ref{prop:intersection equality}, $ T\dot\wedge S=T\wedge S$. Thus it remains to show the first inequality.  
	
	We define a function $\Phi=\Phi(\phi):\CW_1\times \CW_2\to \BR$ as follows:  for every $(\Gamma,\Gamma')\in \CW_1\times \CW_2$, define 
	$$\Phi(\Gamma,\Gamma'):=\sum_{x\in \Gamma\pitchfork \Gamma'} \phi(x).$$
	For every fixed $(\Gamma,\Gamma')$, the transversal intersection points on a compact subset form a finite set.  Since $\phi$ has compact support, the summation in the above definition is a finite summation, hence $\Phi$ is a well-defined non-negative function. 
	
	Next we show that $\Phi$ is lower semi-continuous.  Let $\Gamma_n\to \Gamma$ in $\CW_1$ and $\Gamma'_n\to \Gamma'$ in $\CW_2$.  Let $\left\{x_1,...,x_N\right\}$ be the set of transversal intersection points of $\Gamma$ and $\Gamma'$ in $\supp \phi$.  Then $\Phi(\Gamma,\Gamma')=\sum_{i=1}^N \phi(x_i)$. By the openness of transversal intersection, there exist small open balls  $B_i$ centered at $x_i$ such that $\overline{B_i}\cap \overline{B_j}=\emptyset$ for $i\neq j$, and for all $n$ large enough, $\Gamma_n$ and $\Gamma'_n$ have a unique transversal intersection point $x_{i,n}\in B_i$.   We have $x_{i,n}\to x_i$.  Then for all $n$ large enough, $$\Phi(\Gamma_n,\Gamma'_n)\geq \sum_{i=1}^N\phi(x_{i,n}).$$
	By the continuity of $\phi$,  we have
$$\liminf_{n\to+\infty} \Phi(\Gamma_n,\Gamma'_n)\geq \liminf_{n\to +\infty}\sum_{i=1}^N\phi(x_{i,n})=\sum_{i=1}^N\phi(x_{i})=\Phi(\Gamma,\Gamma').$$
Hence $\Phi$ is lower semi-continuous and non-negative.  Since $\mu_n\otimes\nu_n\to \mu\otimes \nu$ in the weak-* topology, we have
$$\liminf_{n\to+\infty}\int_{\CW_1\times \CW_2} \Phi\; d(\mu_n\otimes \nu_n)\geq \int_{\CW_1\times \CW_2} \Phi\; d(\mu\otimes \nu).$$

Finally, since $$\int_{\CW_1\times \CW_2} \Phi\; d(\mu_n\otimes \nu_n)=\int_X \phi\;d(T_n\dot\wedge S_n)$$ and $$\int_{\CW_1\times \CW_2} \Phi\; d(\mu\otimes \nu)=\int_X \phi\;d(T\dot\wedge S),$$
we get the desired result.
\end{proof}

\medskip

\subsection{Currents that are strongly approximable by algebraic curves}
In this subsection,  let $X$ be a complex smooth projective surface. We fix a K\"ahler form $\omega$ on $X$. Volumes of subvarieties in $X$ are computed with respect to $\omega$.   For any irreducible algebraic curve $C\subseteq X$, we denote by $\genus(C)$  the geometric genus of $C$, and by $\sing(C)$ the set of singular points of $C$.  Let $\pi:\hat{C}\to C$ be the resolution of singularities of $C$.  For $x\in \sing(C)$, let $n_x$ be the number of local irreducible components of $C$ at $x$, i.e. $n_x:=\#\pi^{-1}(x)$. 

\begin{definition}[Condition $*$]\label{defi:condition *}
Let $(C_n)$ be a sequence of irreducible algebraic curves in $X$. We say that this sequence satisfies Condition $*$, if there exists a constant $C>0$ such that for every $n\geq 1$, 
$$\genus(C_n)+\sum_{x\in \sing(C_n)} n_x\leq C\vol(C_n).$$
\end{definition}

\medskip

Let $\Omega\subseteq X$ be an open subset with finitely many connected components. We define uniform laminar currents on $\Omega$ as the sum of uniform laminar currents on the connected components of $\Omega$.

The following theorem is a combination of two results of Dujardin \cite[Theorem 3.1]{DujardinMathann}, \cite[Th\'eor\`eme on page 2]{dujardin04}.  More precisely, assertions~(1) and~(4) are obtained from the construction in \cite[Proposition~4.4]{dujardin04}.
Assertion~(2) follows from the argument of \cite[Proof of Proposition~3.4]{DujardinMathann}. Assertion~(3) is proved in
\cite[Proof of Theorem~4.2]{dujardin04}.

\begin{theorem}\label{thm:dujardin}
Let $(C_n),(D_n)$ be two sequences of irreducible algebraic curves on $X$ satisfying Condition $*$. Assume that \[T_n:=[C_n]/\vol(C_n)\to T\quad \text{and} \quad S_n:=[D_n]/\vol(D_n)\to S\]  in the sense of currents.  Assume that  $T$ and $S$ have continuous local potentials.  Then for every $\ep>0$, by passing to a subsequence,  there exist an open subset $\Omega_\ep\subseteq X$ with finitely many connected components,  two uniformly woven  families $\CW_{1,\ep}$, $\CW_{2,\ep}$ on $\Omega_\ep$,  and finite Borel measures $\mu_\ep,\mu_{n,\ep}$ on $\CW_{1,\ep}$,  $\nu_\ep,\nu_{n,\ep}$ on $\CW_{2,\ep}$, which are all supported on uniformly laminar families,  such that for the corresponding uniformly laminar currents $T_\ep, T_{n,\ep}$ and $S_\ep, S_{n,\ep}$, we have

\smallskip

(1) $T_\ep\leq T,S_\ep\leq S$, and for all $n$, $T_{n,\ep}\leq T_n,S_{n,\ep}\leq S_n$;

\smallskip

(2) $\mu_{n,\ep}\to \mu_\ep$, $\nu_{n,\ep}\to \nu_\ep$ in the weak-* topology;

\smallskip

(3) $\Ma(T_\ep\wedge S_\ep)\geq(1-\ep)\Ma(T\wedge S)$; 

\smallskip

and

\smallskip

(4) There exists an open subset $C_{n,\ep}\subseteq C_n$ (resp. $D_{n,\ep}\subseteq D_n$) contained in the smooth part of $C_n$ (resp. $D_n$) which is closed in $\Omega_\ep$ such that $T_{n,\ep}=[C_{n,\ep}]/\vol(C_n)$ (resp. $S_{n,\ep}=[D_{n,\ep}]/\vol(D_n)$).
\end{theorem}

Although Dujardin’s arguments  \cite[Theorem 3.1]{DujardinMathann} were originally formulated for curves in $\BP^2$,  the same construction applies to curves on an arbitrary smooth projective surface. Indeed, one chooses finitely many generic pencils  $X\dashrightarrow \BP^1$ which provide pairs of transversal local projections, resolves their base points, and applies the good-component estimates chartwise. The required ramification bounds follow from the Riemann–Hurwitz formula (cf. \cite[IV. Corollary~2.4]{Hartshorne}) and Condition $*$.  The geometric intersection argument is local on the surface.  See \cite[Section 1.2]{dujardin23} and \cite[Appendix B]{cantat2023}.

\medskip

The following is the main result of this section.

\begin{theorem}\label{thm:number transversal intersection}
Let $(C_n),(D_n)$ be two sequences of irreducible algebraic curves on $X$ satisfying Condition $*$. Assume that   \[T_n:=[C_n]/\vol(C_n)\to T\quad\text{and}\quad S_n:=[D_n]/\vol(D_n)\to S\] in the sense of currents.  Assume that  $T$ and $S$ have continuous local potentials and $T\wedge S>0$.  Then for all $n$ large enough, we have $\langle C_n,D_n\rangle>0$, and 
$$\lim_{n\to+\infty}\frac{\#(C_n\pitchfork D_n)}{\langle C_n,D_n\rangle}=1,$$
where $C_n\pitchfork D_n$ is the set of transversal intersection points of $C_n\cap D_n$, and $\langle C_n,D_n\rangle$ is the intersection number.
\end{theorem}
\begin{proof}
We first show that for all $n$ large enough, $\langle C_n,D_n\rangle>0$.  Otherwise, by passing to a subsequence, we may assume $\langle C_n, D_n\rangle\leq 0$.  Then in the cohomological level, we have $0\geq \langle\left\{T_n\right\}, \left\{S_n\right\}\rangle\to \langle\left\{T\right\}, \left\{S\right\}\rangle$\footnote{Here and after, $\left\{ \cdot \right\}$ denotes the cohomological class of a closed current (in $\RH^{1,1}(X,\BR)$).}.  But since $T,S$ have continuous local potential, $0<\Ma(T\wedge S)= \langle\left\{T\right\}, \left\{S\right\}\rangle$, a contradiction. Hence $\langle C_n,D_n\rangle>0$ for all $n$ large enough.

 Since $\langle C_n,D_n\rangle>0$, we have $ \#(C_n\pitchfork D_n)\leq \langle C_n,D_n\rangle$. Thus it remains  to show that $$\liminf_{n\to+\infty}\frac{\#(C_n\pitchfork D_n)}{\langle C_n,D_n\rangle}\geq 1.$$
Otherwise, by passing to a subsequence, we may assume that there exists $\ep>0$ such that for all $n\geq 1$, 
\begin{equation}\label{eqn:3.2}
	\frac{\#(C_n\pitchfork D_n)}{\langle C_n,D_n\rangle}\leq 1-3\ep.
\end{equation}

For the fixed $\ep>0$ in \eqref{eqn:3.2}, recall that Theorem \ref{thm:dujardin} gives the following data: an open subset $\Omega_\ep\subseteq X$, two uniformly woven families $\CW_{1,\ep}$, $\CW_{2,\ep}$ on $\Omega_\ep$,  and finite Borel measures $\mu_\ep,\mu_{n,\ep}$ on $\CW_{1,\ep}$,  $\nu_\ep,\nu_{n,\ep}$ on $\CW_{2,\ep}$, which are all supported on uniformly laminar families. Note also that since $T_\ep\leq T$ and $S_\ep\leq S$, $T$ and $S$ having continuous local potential implies that  $T_\ep$ and $S_\ep$ have continuous local potential, see \cite[Lemma 8.2]{Bedford-Lyubich-Smillie}. In particular, they are diffuse and $T_\ep\wedge S_\ep$ is admissible.   By Theorem \ref{thm:dujardin} (3), $$\Ma(T_\ep\wedge S_\ep)\geq(1-\ep)\Ma(T\wedge S).$$ 

Since $\Ma(T\wedge S)<+\infty$, we can choose a non-negative continuous function $\phi:X\to \BR$ with compact support in $\Omega_\ep$ such that  $0\leq \phi \leq 1$ and
\begin{equation*}
\int_X \phi\;d(T_\ep\wedge S_\ep)>(1-2\ep)\Ma(T\wedge S).
\end{equation*}
By Proposition \ref{prop:lower mass estimate}, we have
\begin{equation}\label{eqn:3.3}\begin{split}
    \liminf_{n\to+\infty}\int_X \phi  \;d(T_{n,\ep}\dot\wedge S_{n,\ep})&\geq\int_X \phi  \;d(T_{\ep}\dot\wedge S_{\ep})\\ &=\int_X \phi\;d(T_\ep\wedge S_\ep)>(1-2\ep)\Ma(T\wedge S).
    \end{split}
\end{equation}
By  (\ref{eqn:3.3}) and Theorem \ref{thm:dujardin} (4), we have 
\begin{align}\label{eqn:3.4}
\#(C_n\pitchfork D_n)&\geq  \vol(C_n)\vol(D_n)  \int_X \phi  \;d(T_{n,\ep}\dot\wedge S_{n,\ep})\notag \\
&> (1-2\ep) \vol(C_n)\vol(D_n) \Ma(T\wedge S).
\end{align}
On the other hand,  
$$\Ma(T\wedge S)= \langle\left\{T\right\}, \left\{S\right\}\rangle=\lim_{n\to+\infty} \langle\left\{T_n\right\}, \left\{S_n\right\}\rangle=\lim_{n\to+\infty} \frac{\langle C_n, D_n\rangle}{\vol(C_n) \vol(D_n)}.$$
Combining this with (\ref{eqn:3.4}), we get 
$$	\liminf_{n\to+\infty}\frac{\#(C_n\pitchfork D_n)}{\langle C_n,D_n\rangle}\geq  1-2\ep.$$
This contradicts (\ref{eqn:3.2}), which finishes the proof. 
\end{proof}

		\section{Transversality via lamination and dynamics}\label{sec:trans}
		The goal of this section is to prove Proposition \ref{prop:trans}.  We shall use the theory of laminar currents developed in \S\ref{sec:laminar}, together with the dynamical properties of positive entropy automorphisms on K3 surfaces.
		
		\medskip
		
	In this section, let $X$ be a complex projective K3 surface.  Let $f: X\to X$ be a positive entropy automorphism.  Let $\la_f>1$ be the spectral radius of $f^*$ acting on $H^{1,1}(X,\BR)$.  Fix a K\"ahler form $\omega$ on $X$. Volumes and masses are computed with respect to $\omega$. 
	
	We recall the following results of Cantat \cite{cantat01}. 
	
	\begin{theorem}\label{prop:cantat}
	The following hold:
	
	\smallskip
	
	(1) (cf. \cite[Th\'eor\`eme 2.4, Remarque 2.2]{cantat01}) There exists a unique positive closed current $T^+_f$ (resp. $T^-_f$) of total mass $1$ such that $f^*(T^+_f)=\la_f T^+_f$ (resp. $f_*(T^-_f)=\la_f T^-_f$). 
	
	\smallskip
	
	(2) (cf. \cite[Th\'eor\`eme 3.1 and Section 3.3]{cantat01})  The currents $T^\pm_f$ have continuous local potential, and $T^+_f\wedge T^-_f>0$.
	
	\smallskip
	
	(3) (cf. \cite[Th\'eor\`eme 5.3]{cantat01})  We have $\{T^\pm_f\}\in \operatorname{NS}(X)_{\BR}$. For every irreducible algebraic curve $C\subseteq X$, we have 
	$$\la_f^{-n}[ f^n(C)]\to c^+ T^-_f\quad \text{and}\quad
	\la_f^{-n}[ f^{-n}(C)]\to c^- T^+_f$$
	in the sense of currents, where \[c^\pm:=\frac{\langle \{C\}, \{T^\pm_f\}\rangle}{\langle \{T^+_f\}, \{T^-_f\}\rangle}.\footnote{Here and after we also write $\{\cdot\}$ for the cohomological class of a curve in the Néron--Severi group. Note that \cite{cantat01} assumes $\langle \{T^+_f\}, \{T^-_f\}\rangle=1$}.\]
	
	\smallskip
	
	(4) (cf. \cite[Proposition 6.1]{cantat2020}, \cite[Proof of Lemma 5.1]{cantat01}) For every irreducible algebraic curve $C\subseteq X$, the following are equivalent: \begin{enumerate}
	    \item $\langle \{C\}, \{T^+_f\}\rangle =0$;
        \item $\langle \{C\}, \{T^-_f\}\rangle=0$; \item $C$ is $f$-periodic.
	\end{enumerate}
	\end{theorem}
	\medskip
	
	Combining Theorem \ref{thm:number transversal intersection} and Theorem \ref{prop:cantat}, we get the following result, from which we deduce Proposition \ref{prop:trans}.  
	
	\begin{proposition}
	Let $C, D\subseteq X$ be two non-$f$-periodic irreducible algebraic curves. Let $C_n:=f^n(C)$.  Then $\langle C_n, D\rangle >0$ for all $n$ large enough, and 
	$$\lim_{n\to+\infty}\frac{\#(C_n\pitchfork D)}{\langle C_n,D\rangle}=1.$$
	\end{proposition}

\begin{proof}
By Theorem \ref{prop:cantat} (3) and (4),  
$$\langle\la_f^{-n} C_n, D\rangle \to \langle c^+ \{T_f^-\}, \{D\}\rangle>0.$$
Hence $\langle C_n, D\rangle >0$ for all $n$ large enough.

We shall prove that $$\lim_{n\to+\infty}\frac{\#(C_{2n}\pitchfork D)}{\langle C_{2n},D\rangle}=1.$$ The case where $n$ is odd is similar.  Define $D_n:=f^{-n}(D)$.  Since $f$ is an automorphism, we have 
$$\frac{\#(C_{2n}\pitchfork D)}{\langle C_{2n},D\rangle}=\frac{\#(C_{n}\pitchfork D_n)}{\langle C_{n},D_n\rangle}.$$ Thus we need to show that 
	$$\lim_{n\to+\infty}\frac{\#(C_n\pitchfork D_n)}{\langle C_n,D_n\rangle}=1.$$

By the definition of $C_n$ and $D_n$, these two sequences of curves satisfy Condition $*$ (see Definition \ref{defi:condition *}). By Theorem \ref{prop:cantat} (3),  $$T_n:=[C_n]/\vol(C_n)\to c^+T^-_f,\text{ and }S_n:=[D_n]/\vol(D_n)\to c^-T^+_f,$$ and since $C$ and $D$ are not $f$-periodic, Theorem \ref{prop:cantat} (4) and weak convergence on currents imply that $c^+=c^-=1$. Again by Theorem \ref{prop:cantat} (2), $T^\pm_f$ has continuous local potential, and $T^+_f\wedge T^-_f>0$. By Theorem \ref{thm:number transversal intersection}, 
$$\lim_{n\to+\infty}\frac{\#(C_n\pitchfork D_n)}{\langle C_n,D_n\rangle}=1,$$
which finishes the proof. 
\end{proof}

\medskip
	
	\section{Non-transversality via ramification}\label{sec:non-trans}
	The goal of this section is to prove Proposition \ref{prop:non-trans}. The proof is based on the intersection theory on surfaces, in particular the Hodge Index Theorem. 
	
	\medskip

	As in \S\ref{sec:trans}, let $X$ be a complex projective K3 surface. Let $f: X\to X$ be a positive entropy automorphism.  Let $\la_f>1$ be the spectral radius of $f^*$ acting on $H^{1,1}(X,\BR)$. We now prove the following result, which implies Proposition \ref{prop:non-trans}. 
\begin{proposition}\label{prop:birirred}
  Let $C\subseteq X$ be a non-$f$-periodic irreducible algebraic curve. For $n\in\BZ_{\geq 0}$, let $C_n:=f^n(C)$. Let $Y$ be a smooth projective surface and $\pi:Y\to X$ be a dominant generically finite morphism such that $\deg (\pi)\geq 2$. Assume that there exist an infinite subfamily $(C_n)_{n\in I}$ such that for all $n\in I$, there exists an irreducible component $E_n\subseteq \pi^{-1}(C_n)$ such that $\pi|_{E_n}:E_n\to C_n$ is birational. Let $D\subseteq X$ be a non-$f$-periodic irreducible algebraic curve. Let $Z\subseteq Y$ be an irreducible curve such that $\pi(Z)=D$. Then we have
\[\lim_{\substack{n\to +\infty\\n\in I}}\frac{\langle E_n,Z\rangle}{\langle C_n,D \rangle}=\frac{e}{d}>0,\] where \(d:=\deg(\pi)\) and \(e\) is the degree of \(\pi|_Z:Z\to D\). 
\end{proposition}

\begin{proof}
Throughout the proof, all limits are taken
as \(n\to+\infty\) with \(n\in I\).
By Theorem~4.1, there exists \(c>0\) such that
\begin{equation}\label{eq:C_n}
    \lambda_f^{-n}\{C_n\}\to c\theta_f^-,
\end{equation}
where \(\theta_f^-:=\{T_f^-\}\in\operatorname{NS}(X)_{\mathbb R}\). 
Consequently,
\begin{equation}\label{eq:gamman}
\gamma_n:=\lambda_f^{-n}\pi^*\{C_n\}
\to
\gamma:=c\pi^*\theta_f^-\in\operatorname{NS}(Y)_{\mathbb R}.
\end{equation}
We let
\[
\pi^*C_n=E_n+F_n,
\qquad
\alpha_n:=\lambda_f^{-n}\{E_n\},
\qquad
\beta_n:=\lambda_f^{-n}\{F_n\},
\]
where the divisor \(F_n\) is effective. Thus
\(\alpha_n+\beta_n=\gamma_n\).

Let us fix $A$ an ample class on $Y$. Since the classes \(\alpha_n\) and
\(\beta_n\) are pseudo-effective, in view of \eqref{eq:gamman},
\[
0\leq \max\left(\langle A,\alpha_n\rangle, \langle A,\beta_n\rangle\right)
\leq \langle A,\gamma_n\rangle\leq M,
\]
where \(M>0\) is a uniform constant. The subsets of the
pseudo-effective cone with bounded \(A\)-degree are compact.
Hence both sequences $(\alpha_n)_{n\in I},(\beta_n)_{n\in I}$ are relatively compact.

Let \(\alpha\) be any accumulation point of \((\alpha_n)_{n\in I}\).
Passing to a subsequence, we may assume that
\[
\alpha_n\to\alpha\in\operatorname{NS}(Y)_{\mathbb R}.\]
Thanks to \eqref{eq:gamman}, we have
\[
\beta_n\to\beta:=\gamma-\alpha\in\operatorname{NS}(Y)_{\mathbb R}.\]
We claim that the class \(\alpha\) is nef. In fact, the curves
\(E_n,n\in I\) are pairwise distinct, since the curves \(C_n,n\in I\) are
pairwise distinct. Thus, for every irreducible curve \(B\subset Y\),
\(\langle E_n,B\rangle\geq0\) for all but at most one \(n\).
Passing to the limit gives \(\langle\alpha, B\rangle\geq0\).

The class \(\gamma=c\pi^*\theta_f^-\) is nonzero and nef,
with \(\gamma^2=0\).
Moreover, the class \(\beta\) is pseudo-effective. Therefore,
\[
0\leq\langle\gamma,\alpha\rangle
\leq\langle\gamma,\alpha+\beta\rangle
=\gamma^2=0.
\]
Hence 
\begin{equation}\label{eq:gammaalpha}
\langle\gamma,\alpha\rangle=0.
\end{equation}
By the Hodge index theorem
\cite[Theorem~2.1 and \S2.3.2]{CantatSurvey},
the intersection form on \(\gamma^\perp\) is negative
semidefinite with kernel $\BR\gamma$:  for every $\theta$ that satisfies \(\langle\gamma,\theta\rangle=0\), we have \(\theta^2\leq 0\), and \(\theta^2=0\) is equivalent to $\theta=a\gamma$ for some $a\in \BR$.
In our case, since \(\alpha\) is nef, \(\alpha^2\geq0\). Combining \eqref{eq:gammaalpha}, we get 
\begin{equation}\label{eq:alphagamma2}
    \alpha=a\gamma 
\end{equation}for some \(a\in \BR_{\geq 0}\).

Since \(\pi|_{E_n}\) is birational onto \(C_n\), by \eqref{eq:C_n}, we have
\[
\pi_*\alpha_n=\lambda_f^{-n}\{C_n\}
\to c\theta_f^-.
\]
It follows from \eqref{eq:gamman} and \eqref{eq:alphagamma2} that
\[
c\theta_f^-=\pi_*\alpha=a\pi_*\gamma
=ac\pi_*\pi^*\theta_f^-=acd\theta_f^-.
\]
Thus we have $a=\frac{1}{d}$ and
\[
\alpha=\frac{c}{d}\pi^*\theta_f^-.
\]
The above argument actually works for every convergent subsequence,
and relative compactness thus gives
\begin{equation}\label{eq:alphan} \alpha_n\to\frac{c}{d}\pi^*\theta_f^-.
\end{equation}

Finally, since \(D\) is not \(f\)-periodic,
Theorem~4.1 gives \(\langle\theta_f^-,D\rangle>0\).
By the projection formula and \eqref{eq:alphan}, and since 
\(\pi_*\{Z\}=e\{D\}\),
\[
\begin{aligned}
\lim_{\substack{n\to+\infty\\ n\in I}}
\frac{\langle E_n,Z\rangle}{\langle C_n,D\rangle}
&=
\frac{\frac{c}{d}\,\langle\pi^*\theta_f^-,\{Z\}\rangle}
     {c\,\langle\theta_f^-,\{D\}\rangle}\\
&=
\frac{\langle\theta_f^-,\pi_*\{Z\}\rangle}
     {d\,\langle\theta_f^-,\{D\}\rangle}
=\frac{e}{d}.
\end{aligned}
\]
This finishes the proof.
\end{proof}
	
	\medskip
	
	\section{Examples}\label{se:loxoK3woellipfi}

    In \S\ref{sec:rhoX2noellfib}, we discuss the lattice-theoretic characterization for loxodromic K3 surfaces with $\rho=2,3,4$, which do not admit any elliptic fibration. In \S\ref{se:Wehler} we discuss two explicit classes of K3 surfaces first introduced in \cite{Wehler}, which are concrete examples of K3 surfaces over the ground field admitting positive entropy automorphisms. 
    
    \subsection{Loxodromic K3 surfaces without elliptic fibrations}\label{sec:rhoX2noellfib}
    
For convenience, we say that a lattice is \emph{admissible} if it can be realized as the N\'eron--Severi lattice of a (complex) K3 surface. Furthermore, we say that an admissible lattice is \emph{loxodromic} if it corresponds to a (complex) loxodromic K3 surface. Note that by \cite[Remark 17.2.16]{Huybrechts}, every admissible lattice is realizable by a K3 surface over a certain number field.

According to  \cite[\S7]{PyatetskiiShapiroShafarevich} (cf. also \cite[Example 15.2.11]{Huybrechts} or \cite[Theorem 1.2 (2)]{Yu}), a  K3 surface $X$ with $\rho(X)=2$ is loxodromic if and only if there exists no non-zero class $v\in \operatorname{NS}(X)$ such that $v^2=0$ or $-2$. In particular, by \cite[Proposition 11.1.3]{Huybrechts}, such a surface does not admit any elliptic fibration. Let
\[
\begin{pmatrix}2a&b\\ b&2c\end{pmatrix}
\] be an admissible lattice. Then it is loxodromic if and only if the discriminant $b^{2}-4ac$ is not a square, and the quadratic equation $ax^{2}+bxy+cy^{2}=-1$ has no integral solution in $(x,y)$. For instance, the Cayley--Oguiso lattice (cf. \cite{Oguiso,FestiGarbagnatiVanGeemenVanLuijk}) \[
\begin{pmatrix}
4&2\\
2&-4
\end{pmatrix}\] is the N\'eron--Severi lattice of a loxodromic K3 surface, admitting an automorphism $f$ with $\lambda_f=9+4\sqrt{5}$, but with no elliptic fibration.
  
Loxodromic admissible lattices with rank $\geq 3$ are completely classified by Yu \cite[Theorem 1.2 (1)]{Yu}: An admissible lattice of rank $\geq 3$ is loxodromic if and only if it is not isometric to any
of the \(311\) exceptional lattices in that theorem.
Among these loxodromic lattices, the anisotropic ones,
i.e.\ those containing no nonzero vector of square zero,
correspond precisely to K3 surfaces without genus-one
fibrations.

\subsection{Wehler surfaces}\label{se:Wehler}
	By a Wehler-type surface, we mean either the complete intersection defined by a $(1,1)$-form with a $(2,2)$-form in $\BP^2\times\BP^2$, which we shall call of \emph{Type I}, or defined by a single $(2,2,2)$-form in $(\mathbb{P}^1)^3$, which we shall call of \emph{Type II}. 
    
    \subsubsection{Type I} Silverman \cite{Silverman} pioneered the study of rational points via dynamics on these surfaces. We refer to his work for more details (cf. also  \cite[Example 3.5]{Kawaguchi}). 
    
    Type I Wehler surfaces form an $18$-dimensional family. A general member has Picard rank two, and the intersection matrix of the Néron--Severi lattice is (cf. \cite[15.2.5 (i)]{Huybrechts}) $$\left(\begin{matrix}
	2 & 4 \\
	4 & 2 
	\end{matrix}\right).$$ It can be checked that (see \S\ref{sec:rhoX2noellfib}) there is no class of square zero. Hence, such a generic type I Wehler surface does not admit an elliptic fibration.
	
    \medskip

	Let $S$ be a smooth type I Wehler surface. Let $p_i:S\to \mathbb{P}^2,i=1,2$ be the projections onto the two factors of $\mathbb{P}^2\times\mathbb{P}^2$. Assume that both projections $p_i \colon S \to \mathbb{P}^2$, $i=1,2$, are finite (of degree two).
These projections induce involutions $\sigma_i: S\to S$ (that is, an automorphism of order two). The subgroup $\mathcal{A}\subset\operatorname{Aut}(S)$ generated by $\sigma_1,\sigma_2$ is isomorphic to $\BZ/2 *\BZ/2$.
	Let $D_i:=p_i^*\CO(1)\in \operatorname{NS}(S),i=1,2$. Define $\Lambda:=\mathbb{Z}D_1\oplus \mathbb{Z}D_2$, which is a sublattice of $\operatorname{NS}(S)$.  In terms of the basis $\mathbf{e}:=(D_1,D_2)^t$, the action of $\mathcal{A}$ on $\Lambda$ is given by
	$$\sigma_1^*\mathbf{e}=\left(\begin{matrix}
	1 & 0 \\
	4 & -1 
	\end{matrix}\right)\mathbf{e},\quad \sigma_2^* \mathbf{e}=\left(\begin{matrix}
	-1 & 4 \\
	0 & 1 
	\end{matrix}\right)\mathbf{e}.$$
	Let $$\sigma:=\sigma_2\circ\sigma
    _1\in\operatorname{Aut}(S)$$ and $$E_1:=(2+\sqrt{3}) D_1-D_2,\quad E_2:= -D_1+(2+\sqrt{3})D_2$$ as elements of $\operatorname{NS}(S)_\mathbb{R}$. Then $$\sigma^*(E_1)=(7+4\sqrt{3})E_1,\quad \sigma^*(E_2)=(7+4\sqrt{3})^{-1}E_2.$$ Hence the automorphism $\sigma$ has positive entropy, and by Theorem \ref{thm:main} such an $S$ has the potential Hilbert property.

The class $E_1+E_2$ is ample, and hence by \cite[Proposition 3.1 (2)]{Kawaguchi}, there exist no $\sigma$-periodic curves on $S$. Therefore, if $S$ contains a rational curve birational to $\BP^1$ over the ground field, then Theorem \ref{thm:hilbbfempty} implies that $S$ has the Hilbert property (over the ground field). 

For an explicit numerical example, let $S\subset\BP_{\BQ}^2\times\BP_{\BQ}^2$ be defined by the common zeroes of the $(1,1)$-form \[
x_0y_0+x_1y_1+x_2y_2,
\]
and the $(2,2)$-form
\[
\begin{aligned}
&-2x_0^2y_0^2-x_0^2y_1y_2
-2x_0x_1y_0y_1-3x_0x_2y_0y_2-x_0x_2y_1^2\\
&+2x_1^2y_2^2+2x_1x_2y_0^2
-2x_1x_2y_0y_1-x_1x_2y_1^2\\
&-x_2^2y_0y_1-x_2^2y_0y_2-x_2^2y_1y_2.
\end{aligned}
\]
Then it can be checked that $S$ is smooth, the projections onto the two factors $\BP^2$ are finite, and $S$ contains the rational curve
\[
x_0y_0+x_1y_1=0,\quad x_2=0,\quad y_2=0.
\]
	
    \subsubsection{Type II} The dynamics of such surfaces were studied in depth by McMullen \cite{McMullen} (cf. also \cite[\S1.2]{CantatSurvey}) and the relation with rational points was addressed by the works \cite{Baragar,Billard,Wang}.

    Type II Wehler surfaces form an $18$-dimensional family. A general member has Picard rank $3$ and the intersection matrix of the Néron--Severi lattice is
    $$\left(\begin{matrix}
	0 & 2 & 2 \\
	2 & 0 & 2 \\
    2 & 2 & 0
	\end{matrix}\right).$$
	
	\medskip

	Let $S$ be a smooth type II Wehler surface.
    For $i=1,2,3$, let $p_i:S\to\BP^1$ be induced by the projection $(\BP^1)^3\to\BP^1$ onto the $i$-th factor. The general fibre of $p_i:S\to\BP^1$ is defined by a $(2,2)$-form in $(\BP^1)^2$, thus a genus one curve. Hence, the potential Hilbert property of $S$ follows from \cite{GCM24}.

For $j=1,2,3$, let $\pi_{j}:S\to \mathbb{P}^1\times \BP^1$ be induced by the projection $(\BP^1)^3\to (\BP^1)^2$ without the $j$-th factor. Assume that every projection
$\pi_j,j=1,2,3$, is finite (of degree two). It therefore induces an involution $\sigma_j:S\to S$ defined over the ground field. The subgroup $\mathcal{A}\subset\operatorname{Aut}(S)$ generated by $\sigma_1,\sigma_2,\sigma_3$ is isomorphic to  $\BZ/2 *\BZ/2*\BZ/2$. For $i=1,2,3$, let $D_i=p_i^*\{\infty\}\in \operatorname{NS}(S)$.
	In terms of the basis $\mathbf{e}=(D_1,D_2,D_3)^t$, the action of $\CA$ admits the matrix representation:
	$$\sigma_1^* \mathbf{e}=\left(\begin{matrix}
		-1 & 2 & 2\\
		0 & 1 & 0 \\
		0 & 0 & 1 
	\end{matrix}\right)\mathbf{e},\quad \sigma_2^* \mathbf{e}=\left(\begin{matrix}
		1 & 0 &0\\
		2 & -1 & 2 \\
		0 & 0 & 1 
	\end{matrix}\right)\mathbf{e},\quad\sigma_3^*\mathbf{e}=\left(\begin{matrix}
		1 & 0 & 0\\
		0 & 1 & 0\\
		2 & 2 & -1 
	\end{matrix}\right)\mathbf{e} .$$

	Let $$E_1=\frac{1-\sqrt{5}}{2}D_1+\frac{1+\sqrt{5}}{2}D_2+D_3;$$
	$$E_2=\frac{1+\sqrt{5}}{2}D_1+\frac{1-\sqrt{5}}{2}D_2+D_3$$
    as elements of $\operatorname{NS}(S)_\BR$.
	Let $$\sigma:=\sigma_2\circ\sigma_3\circ \sigma_1\in\operatorname{Aut}(S).$$
	Then we have 
	$$\sigma^*(E_1)=(9+4\sqrt{5})E_1,\quad \sigma^*(E_2)=(9+4\sqrt{5})^{-1}E_2.$$
	Hence the automorphism $\sigma$ has positive entropy, and the potential Hilbert property can also be deduced from Theorem \ref{thm:main}.
    
	Since the class $E_1+E_2$ is ample, it follows from \cite[Proposition 3.1 (2)]{Kawaguchi} that there are no $\sigma$-periodic curves on $S$. 
	 Note that each genus one fibration $p_i$ necessarily has singular fibres. If one of them has an irreducible component which is birational to $\BP^1$ over the ground field, then Theorem \ref{thm:hilbbfempty} implies that $S$ has the Hilbert property (over the ground field).

     For a numerical example, consider the affine model
\[
x^2+y^2+z^2+x^2y^2+y^2z^2+z^2x^2
-3x^2y^2z^2-4xyz=1
\]
in $(\mathbb A_{\mathbb Q}^1)^3$. Its closure in $(\mathbb P_{\mathbb Q}^1)^3$
is smooth, and all three projections onto two factors are finite
of degree two. The fibre at $z=1$ has equation
\[
x^2y^2=(x-y)^2.
\]
Its two irreducible components are given by
\[
xy=x-y
\qquad\text{and}\qquad
xy=y-x,
\]
and both are rational curves birational to $\mathbb P^1_{\mathbb Q}$.

    \medskip
    
	 \section*{Acknowledgments} 
	The collaboration was initiated during the workshop ``Frontiers in Arithmetic Dynamics'' (FIAD2024) held at Beijing International Center for Mathematical Research. Substantial progress was made during a visit to Westlake University in the summer of 2024 and during the 2025 Summer Research Institute in Algebraic Geometry at Colorado State University. The hospitality is greatly appreciated.  
    
    This work was reported by Z. H. at ``Conference on Algebraic Geometry'' held at Tsinghua Sanya International Forum in December 2025, and by J. X. at Harbin Institute of Technology (HIT) in January 2026. We acknowledge the organizers for their kind invitation. Z. H. thanks Pietro Corvaja for interesting discussions.

    \subsection*{AI disclosure} Figure \ref{fig:xtox2} and the numerical examples in \S\ref{se:Wehler} were suggested by ChatGPT 5.6. The authors independently verified the correctness of these examples. All other statements and proofs in the paper were developed by the authors, who take full responsibility for the contents of the paper.
    
\subsection*{Funding} Damián Gvirtz-Chen is supported by UKRI award UKRI094. Zhizhong Huang is supported by National Key R\&D Program of China (No. 2025YFA1017300 and 2025YFA1017303). Zhuchao Ji is supported by the National Key R\&D Program of China (No. 2025YFA1018300), NSFC Grant (No.12401106),
and ZPNSF Grant (No. XHD24A0201). Junyi Xie is supported by NSFC Grant (No. 12271007)  and by the Xplorer Prize from the New Cornerstone Science Foundation.
	
	\bibliographystyle{alpha}
	\bibliography{posentropy}
	
\end{document}